\input amstex
\documentstyle{amsppt}
\pagewidth{6.5in}
\topmatter
\title
FASTER AND FASTER CONVERGENT SERIES FOR $\zeta(3)$
\endtitle
\author Tewodros Amdeberhan
\endauthor
\affil Department of Mathematics, Temple University, 
Philadelphia PA 19122, USA \\
tewodros\@euclid.math.temple.edu
\endaffil
\date Submitted: April 8, 1996. Accepted: April 15, 1996
\enddate
\abstract 
Using WZ pairs we present accelerated series for computing $\zeta(3)$
\endabstract
\endtopmatter
\def\({\left(}
\def\){\right)}

\font\smcp=cmcsc8
\headline={\ifnum\pageno>1{\smcp the electronic journal of combinatorics 3 (1996), \#R13\hfill\folio} \fi}
\midspace{.1in}

Alf van der Poorten [P] gave a delightful account of Ap\'ery's proof [A] of the irrationality of $\zeta(3)$. Using WZ forms, that came from [WZ1], Doron Zeilberger [Z] embedded it in a conceptual framework.\smallskip

We recall [Z] that a discrete function A(n,k) is called Hypergeometric (or Closed Form (CF)) in two variables  when the ratios $A(n+1,k)/A(n,k)$ and $A(n,k+1)/A(n,k)$ are both rational functions. A pair (F,G) of CF functions is a WZ pair if $F(n+1,k) - F(n,k) = G(n,k+1) - G(n,k)$. In this paper, after choosing a particular F (where its companion G is then produced by the amazing Maple package EKHAD accompanying [PWZ]), we will give a list of accelerated series calculating $\zeta(3)$. Our choice of F is
$$F(n,k) = \frac {(-1)^k k!^2 (sn-k-1)!}{(sn+k+1)!(k+1)}$$
where s may take the values s=1,2,3, \dots [AZ] (the section pertaining to this can be found in \linebreak http://www.math.temple.edu/\~{}tewodros). In order to arrive at the desired series we apply the following result:\smallskip
\bf Theorem: \rm ([Z], Theorem 7, p.596) For any WZ pair (F,G)
$$\sum_{n=0}^{\infty}G(n,0) = \sum_{n=1}^{\infty}\(F(n,n-1)+G(n-1,n-1)\),$$
whenever either side converges.\smallskip
The case s=1 is Ap\'ery's celeberated sum [P] (see also [Z]):
$$\zeta(3) = \frac 52 \sum_{n=1}^{\infty}(-1)^{n-1}\frac 1{\binom {2n}n n^3}$$
where the corresponding G is
$$G(n,k)=\frac {2(-1)^k k!^2 (n-k)!}{(n+k+1)!(n+1)^2}.$$

\pagebreak
\midspace{.5in}

For s=2 we obtain 
$$
\align
\zeta(3) &=
\frac 14 \sum_{n=1}^{\infty}(-1)^{n-1} \frac {56n^2-32n+5}{(2n-1)^2} \frac 1{\binom {3n}n \binom {2n}n n^3}
\endalign$$
where G is
$$G(n,k)= \frac {(-1)^k k!^2 (2n-k)!(3+4n)(4n^2+6n+k+3)}{2(2n+k+2)!(n+1)^2(2n+1)^2}.$$

For s=3 we have 
$$\zeta(3) =
\sum_{n=0}^{\infty}\frac {(-1)^n}{72\binom {4n}n \binom {3n}n }\{\frac{6120n+5265n^4+13761n^2+13878n^3+1040}{(4n+1)(4n+3)(n+1)(3n+1)^2(3n+2)^2}\},
$$
and so on.

\enddocument